Scientific
Research
Publishing

# A Class of Smoothing Modulus-Based Iterative Method for Solving Implicit Complementarity Problems

## Cong Guo, Chenliang Li, Tao Luo

School of Mathematics and Computing Science, Guilin University of Electronic Technology, Guilin, China
Email: chenli@guet.edu.cn

## Abstract

In this paper, a class of smoothing modulus-based iterative method was presented for solving implicit complementarity problems. The main idea was to transform the implicit complementarity problem into an equivalent implicit fixed-point equation, then introduces a smoothing function to obtain its approximation solutions. The convergence analysis of the algorithm was given, and the efficiency of the algorithms was verified by numerical experiments.

## 1. Introduction

Complementarity problems arise widely in many applications of science and engineering, such as elastic contact, economic transport, boundary problems in fluid dynamics, convex quadratic programming and inverse problems [1]. In this paper, we consider the implicit complementarity problems (ICP) [2]: find vectors $z, w \in R^n$, such that

$$z - m(z) \geq 0, \quad w = Mz + q \geq 0$$

$$(z - m(z))^\mathrm{T} (Mz + q) = 0 .$$

where $M \in R^{n \times n}$ is a given matrix, $q \in R^n$ is a given constant vector, $m(\cdot)$ is a mapping from $R^n$ into itself, and if the mapping $m(\cdot)$ is a zero mapping, the ICP reduces to linear complementarity problem (LCP).

ICP was introduced into the complementarity theory as a mathematical tool in the study of some stochastic optimal control problems, it was first proposed by Bensoussan and Lion in [3]. After decades of research and development, people have obtained the existence complementarity theorem of solutions to the implicit complementarity problem [4], and proposed many effective methods for solving implicit complementarity problems. Such as fixed-point method [5] [6], projection iteration method [7] [8], Schwarz method [9] etc. Apart from that Zheng and Qu [10] established a hybrid method for solving implicit complementarity problems with superconvergence, numerical examples show that this method has higher accuracy and faster convergence than some existing methods. Tian *et al.* [11] proposed an unconstrained and differentiable penalty method for solving the implicit complementarity problems. Xia and Li [12] firstly constructed modulus-based matrix splitting iterative method for solving nonlinear complementarity problems. Then Hong and Li [13] presented modulus-based matrix splitting iterative method (MMS) for solving the implicit complementarity problems, and analyzes its convergence. Wang *et al.* [14] given new convergence conditions of MMS when the system matrix is a positive-definite matrix and an $H_+$-matrix. On their basis, Yin *et al.* [15] proposed a class of accelerated modulus-based matrix splitting iteration methods to solve the implicit complementarity problems. Zheng and Vong [16] proposed a modified modulus-based matrix splitting iterative method to solve the implicit complementarity problems. Wang and Cao [17] constructed a two-step modulus-based matrix splitting iterative methods (TMMS) for solving implicit complementarity problems, numerical experiments show that this iterative method is more effective than the MMS methods.

Dong and Jiang proposed modular iteration method in [18]. Its main idea is to transform the linear complementarity problems into an equivalent fixed-point equation system, and then solve it. Foutayeni *et al.* [19] constructed a smoothing function to approximate the linear complementarity problems, and proposed a class of modified Newton methods and $m+1$-step iteration methods to solve the linear complementarity problems, obtaining an effective smoothing numerical algorithm. In this paper, we extend this method to solve the implicit complementarity problems. Firstly, we transform the implicit complementarity problems into an equivalent implicit fixed-point equation system. Since it is a non-differentiable system of absolute value equations, we introduce a smooth function to approximate the original problem. Then we construct a class of smoothing modulus-based iteration method for solving the approximated system of equations. Finally, we analyze the convergence of the new method, and verify its effectiveness by numerical experiments.

The organization of the paper is as follows. In Section 2, we establish the smoothing modulus-based iteration method for solving the implicit complementarity problems. The convergence of the smoothing modulus-based iteration method is presented in Section 3, and the numerical results about the new methods are shown and discussed in Section 4. In addition, some conclusion is given in Section 5.

## 2. The Smoothing Modulus-Based Iterative Method

In this section, we consider the following implicit complementarity problems (ICP): Given matrix $M \in R^{n \times n}$ and constant vector $q \in R^n$, and $m(\cdot)$ is a mapping from $R^n$ into itself, find vectors $z, w \in R^n$, such that

$$z - m(z) \geq 0, \quad w = Mz + q \geq 0 \tag{1}$$
$$(z - m(z))^{\mathrm{T}} (Mz + q) = 0.$$

Let $z - m(z) = \beta(|x| + x)$, $w = \alpha(|x| - x)$, transform the implicit complementarity problems into a fixed-point equation:

$$(\alpha I + \beta M) x = (\alpha I - \beta M)|x| - M \cdot m(z) - q, \tag{2}$$

where $\alpha, \beta$ are two positive constants and $I \in R^{n \times n}$ is the identity matrix.

Let $g(z) = z - m(z)$. Since $g(z)$ is an invertible function, then $z = g^{-1}(\beta(|x| + x))$. The implicit fixed-point equation of (2) is changed as following,

$$(\Omega + A) x = (\Omega - A)|x| - M \cdot m(g^{-1}(\beta(|x| + x))) - q, \tag{3}$$

where $\Omega = \alpha I \in R^{n \times n}$, $A = \beta M \in R^{n \times n}$.

**Lemma 1.** (a) If the solution of the implicit complementarity problems (1) is $(z, w)$, then $x = \frac{1}{2}\left(\dfrac{z - m(z)}{\beta} - \dfrac{w}{\alpha}\right)$ is the solution of (3).

(b) If the vector $x$ satisfies (3), then $z - m(z) = \beta(|x| + x)$, $w = \alpha(|x| - x)$ is the solution of (1).

**Proof.** According to Theorem 2.1 in [13], we can similarly prove it.

For the implicit fixed-point Equation (3), let $F(x)$ be a vector function:

$$F(x) = (\Omega + A) x - (\Omega - A)|x| + M \cdot m(z) + q, \tag{4}$$

where $\Omega = \alpha I \in R^{n \times n}$, $A = \beta M \in R^{n \times n}$.

Due to $|x|$, $F(x)$ is non-differentiable. We introduce a smoothing vector function from $R^n \to R^n$ [19]

$$\left(x^2 + \mathrm{e}^{-c}\right)^{\frac{1}{2}} = \left(\left(x_1^2 + \mathrm{e}^{-c}\right)^{\frac{1}{2}}, \left(x_2^2 + \mathrm{e}^{-c}\right)^{\frac{1}{2}}, \cdots, \left(x_n^2 + \mathrm{e}^{-c}\right)^{\frac{1}{2}}\right)^{\mathrm{T}},$$

where $c$ is a large positive integer.

Substitute it into the function (4), we can get a smoothing function:

$$F_c(x) = (\Omega + A) x - (\Omega - A)\left(x^2 + \mathrm{e}^{-c}\right)^{\frac{1}{2}} + M \cdot m(z) + q,$$

where $\Omega = \alpha I \in R^{n \times n}$, $A = \beta M \in R^{n \times n}$.

**Lemma 2.** When $c \to \infty$, $F_c(x)$ uniformly converges to $F(x)$.

**Proof.** According to Corollary 2.1 in [19], we similarly prove it. When $c \to \infty$, $\left(x^2 + \mathrm{e}^{-c}\right)^{\frac{1}{2}}$ converges uniformly to $|x|$. Therefore, $F_c(x)$ uniformly converges to $F(x)$.

From Lemma 2, we know that if $x^*$ is the solution of $F_c(x) = 0$, then $x^*$

converges to the solution of $F(x) = 0$. And $F_c(x) = 0$ is a smooth nonlinear system of equations, so we can use classical Newton method to solve it, in order to get a better initial value, we use following modulus-based iteration method.

**Algorithm 1.** (Modulus-Based Iteration Method)

Step 1: Given $\varepsilon > 0$, $z^0 \in R^n$, set $k = 0$.

Step 2: 1) Calculate the initial vector $x^0 = \frac{1}{2}\left(z^k - w^k - m(z^k)\right)$, set $j = 0$,

2) Iterative Computing $x^{j+1} \in R^n$ by solving the equations

$$(\Omega + A)x^{j+1} = (\Omega - A)\left|x^j\right| - M \cdot m(z^k) - q \tag{5}$$

3) $z^{k+1} = \beta\left(\left|x^{j+1}\right| + x^{j+1}\right) + m(z^k)$

Step 3: If $\text{RES}(z^k) = \left\|\left(Mz^k + q\right)^T\left(z^k - m(z^k)\right)\right\| < \varepsilon$, then stop. Otherwise, set $k = k + 1$ and return to Step 2.

According to Remark 1 in [13], $x^k \to x^*$, when $k \to \infty$. Therefore, there exists some $k_0$, $\left\|x^{k_0} - x^*\right\| < \delta$, $\delta > 0$ is a constant. That means that $x^{k_0}$ is a better approximation vector. So we can construct the following new algorithm.

**Algorithm 2.** (Smoothing Modulus-Based Newton Method)

Step 1: Given $\varepsilon > 0$, $c > 0$, set $k = 0$.

Step 2: Get the initial vector $x^0$ through Algorithm 1.

Step 3: Compute $F_c(x^k)$,

$$F_c'(x^k) = \begin{pmatrix} \dfrac{\partial f_1(x^k)}{\partial x_1} & \dfrac{\partial f_1(x^k)}{\partial x_2} & & \dfrac{\partial f_1(x^k)}{\partial x_n} \\ \dfrac{\partial f_2(x^k)}{\partial x_1} & \dfrac{\partial f_2(x^k)}{\partial x_2} & \cdots & \dfrac{\partial f_2(x^k)}{\partial x_n} \\ \vdots & \vdots & \ddots & \vdots \\ \dfrac{\partial f_n(x^k)}{\partial x_1} & \dfrac{\partial f_n(x^k)}{\partial x_2} & \cdots & \dfrac{\partial f_n(x^k)}{\partial x_n} \end{pmatrix},$$

$$\Delta x^k = -\left(F_c'(x^k)\right)^{-1} F_c(x^k).$$

Step 4: Compute

$$x^{k+1} = x^k + \Delta x^k, \quad z^{k+1} = \beta\left(\left|x^{k+1}\right| + x^{k+1}\right) + m(z^k).$$

Step 5: If $\left\|x^{k+1} - x^k\right\|_2 < \varepsilon$, then stop and output

$$z = \beta\left(\left|x^{k+1}\right| + x^{k+1}\right) + m(z^{k+1}),$$

$$w = \alpha\left(\left|x^{k+1}\right| - x^{k+1}\right).$$

Otherwise, let $x^k = x^{k+1}$, $z^k = z^{k+1}$, $k = k + 1$, return to Step 2.

The classic Newton method needs to choose a good initial vector for better convergence. We introduce the parameter $\xi_1, \xi_2$, and use the following iterative sequence to improve the classical Newton method.

$$\begin{cases} y^k = x^k \pm \xi_1 F_c'(x^k)^{-1} F_c(x^k), \\ x^{k+1} = y^k - \xi_2 F_c'(x^k)^{-1} F_c(y^k). \end{cases} \tag{6}$$

We take $\xi_1 = \xi_2 = 1$, and use the modified Newton method to solve $F_c(x) = 0$.

**Algorithm 3.** (Modified Smoothing Modulus-based Newton method)

Step 1: Given $\varepsilon > 0$, $c > 0$, set $k = 0$.

Step 2: Get the initial vector $x^0$ through Algorithm 1.

Step 3: Computing $F_c(x^k)$, $F'_c(x^k)$.

Step 4: Use (6) compute $x^{k+1}$, and

$$z^{k+1} = \beta\left(\left|x^{k+1}\right| + x^{k+1}\right) + m\left(z^k\right).$$

Step 5: If $\left\|x^{k+1} - x^k\right\|_2 < \varepsilon$, then stop and output

$$z = \beta\left(\left|x^{k+1}\right| + x^{k+1}\right) + m\left(z^{k+1}\right),$$

$$w = \alpha\left(\left|x^{k+1}\right| - x^{k+1}\right).$$

Otherwise, let $x^k = x^{k+1}$, $z^k = z^{k+1}$, $k = k+1$, return to Step 2.

On the basis of the modified Newton method, we do $m+1$-step iterations in the iterative process, and then construct a smoothing modulus-based $m+1$-step iteration method to solve the implicit complementarity problem, and the iterative sequence is as follows:

$$\begin{cases} y_1^k = x^k \pm F'_c\left(x^k\right)^{-1} F\left(x^k\right) \\ y_2^k = y_1^k - F'_c\left(x^k\right)^{-1} F\left(y_1^k\right) \\ \qquad\qquad \vdots \\ y_m^k = y_{m-1}^k - F'_c\left(x^k\right)^{-1} F\left(y_{m-1}^k\right) \\ x^{k+1} = y_m^k - F'_c\left(x^k\right)^{-1} F\left(y_m^k\right) \end{cases} \tag{7}$$

we use the $m+1$-step iterative method to solve $F_c(x) = 0$.

**Algorithm 4.** (Smoothing Modulus-Based $m+1$-step Iterative Method)

Step 1: Given $\varepsilon > 0$, $c > 0$, set $k = 0$.

Step 2: Get the initial vector $x^0$ through Algorithm 1.

Step 3: Computing $F_c(x^k)$, $F'_c(x^k)$.

Step 4: Use (7) compute $x^{k+1}$, and

$$z^{k+1} = \beta\left(\left|x^{k+1}\right| + x^{k+1}\right) + m\left(z^k\right).$$

Step 5: If $\left\|x^{k+1} - x^k\right\|_2 < \varepsilon$, then stop and output

$$z = \beta\left(\left|x^{k+1}\right| + x^{k+1}\right) + m\left(z^{k+1}\right),$$

$$w = \alpha\left(\left|x^{k+1}\right| - x^{k+1}\right).$$

Otherwise, let $x^k = x^{k+1}$, $z^k = z^{k+1}$, $k = k+1$, return to Step 2.

## 3. Convergence Theorem

In this section, we give the convergence of the above algorithms.

**Theorem 1.** If $x^*$ is the solution of the system of equations $F_c(x) = 0$, the iterative sequence $\left\{x_c^k\right\}$ is second-order convergent in a neighborhood of $x_c^*$.

**Proof.** According to the convergence theorem of Newton method, it can be obtained that the iterative sequence generated by Algorithm 2 converges to the solution $x^*$ of the equation $F_c(x) = 0$ with the order two in a neighborhood of $x_c^*$.

**Theorem 2.** For nonlinear vector functions $F_c(x) = 0$, define the sequence:

$$\begin{cases} y^k = x^k \pm \xi_1 F_c'(x^k)^{-1} F_c(x^k), \\ x^{k+1} = y^k - \xi_2 F_c'(x^k)^{-1} F_c(y^k) \end{cases}$$

if $x^*$ is the solution of the system of equations $F_c(x) = 0$, when $\xi_1 = \pm 1$, $\xi_2 = 1$, Algorithm 3 converges to $x^*$ with order three.

**Proof** The Taylor expansion of the function $F_c(x^k)$ at $x^*$:

$$F_c(x^k) = F_c'(x^*)(x^k - x^*) + \frac{1}{2} F_c''(x^*)(x^k - x^*)^2 + o\left(\left\| x^k - x^* \right\|^3\right). \tag{8}$$

For (8), let $r^k = x^k - x^*$, we have

$$F_c(x^k) = F_c'(x^*)\left(r^k + \frac{1}{2} F_c''(x^*)(r^k)^2 + o\left(\left\| r^k \right\|^3\right)\right), \tag{9}$$

$$F_c'(x^k) = A = F_c'(x^*)\left[I + F_c'(x^*)^{-1} F_c''(x^*) r^k + o\left(\left(r^k\right)^2\right)\right]. \tag{10}$$

Then

$$A^{-1} = \left[I - F_c'(x^*)^{-1} F_c''(x^*) r^k + o\left(\left(r^k\right)^2\right)\right] F_c'(x^*)^{-1}. \tag{11}$$

Let $r_y^k = y^k - x^*$, According to (11), that

$$r_y^k = y^k - x^* - \xi_1 A^{-1} F_c'(x^k) = (1 - \xi_1) r^k + \frac{\xi_1}{2} F_c'(x^*)^{-1} F_c''(x^*)(r^k)^2 + o\left(\left(r^k\right)^3\right). \tag{12}$$

Which leads to

$$\begin{aligned} A r_y^k = F_c'(x^*) \cdot &\left[I + F_c'(x^*)^{-1} F_c''(x^*) r^k + o\left(\left(r^k\right)^2\right)\right](1 - \xi_1) r^k \\ &+ \frac{\xi_1}{2} F_c'(x^*)^{-1} F_c''(x^*)(r^k)^2 + o\left(\left(r^k\right)^3\right). \end{aligned} \tag{13}$$

From (8),

$$F_c(y^k) = F_c'(x^*) \cdot \left(r_y^k + \frac{1}{2} F_c'(x^*)^{-1} F_c''(x^*)(r_y^k)^2 + o\left(\left(r_y^k\right)^3\right)\right). \tag{14}$$

On the other hand

$$r_y^{k+1} = x^{k+1} - x^* = y^k - x^* - \xi_2 A^{-1} F_c(y^k). \tag{15}$$

Equivalent to

$$A r_y^{k+1} = A r_y^k - \xi_2 F_c(y^k). \tag{16}$$

Applying formula (11), (13) and (14), we have

$$A r_y^{k+1} = F_c'(x^*) \cdot \left(I + F_c'(x^*)^{-1} F_c''(x^*) r^k + o\left(\left(r^k\right)^2\right)\right) r^{k+1}, \tag{17}$$

$$Ar_y^k - \xi_2 F_c\left(y^k\right) = F_c'\left(x^*\right) \cdot \left(\left(1-\xi_1\right)\left(1-\xi_2\right)r_y^k + \left(1 - \frac{\xi_1}{2} - \xi_2\left(\frac{1-\xi_1+\xi_1^2}{2}\right)\right)\right)$$
$$\cdot F_c'\left(x^*\right)^{-1} F_c''\left(x^*\right)\left(r^k\right)^2 + o\left(\left(r^k\right)^3\right)\right). \tag{18}$$

In order to get $r^{k+1} = o\left(\left(r^k\right)^3\right)$, the constants $\xi_1$ and $\xi_2$ must satisfy:

$$\begin{cases} \left(1-\xi_1\right)\left(1-\xi_2\right) = 0 \\ 1 - \frac{\xi_1}{2} - \xi_2\left(\frac{1-\xi_1+\xi_1^2}{2}\right) = 0 \end{cases}, \tag{19}$$

by solving the above equations, we can get $\xi_1 = \pm 1, \xi_2 = 1$. The theorem is proved.

**Theorem 3.** If $x^*$ is the solution of the system of equations $F_c\left(x\right) = 0$, Sequence

$$\begin{cases} y_1^k = x^k \pm F_c'\left(x^k\right)^{-1} F\left(x^k\right) \\ y_2^k = y_1^k - F_c'\left(x^k\right)^{-1} F\left(y_1^k\right) \\ \quad\quad \vdots \\ y_m^k = y_{m-1}^k - F_c'\left(x^k\right)^{-1} F\left(y_{m-1}^k\right) \\ x^{k+1} = y_m^k - F_c'\left(x^k\right)^{-1} F\left(y_m^k\right) \end{cases} \tag{20}$$

for any positive integer *m*, the sequence produced by Algorithm 4 is converges to $x^*$ with order $m+2$.

**Proof.** We use the mathematical induction method to prove it. Obviously, when $m=1$, the sequence is third-order convergent from Theorem 2, and the theorem holds. Then we assume that $m-1 \geq 1$ holds, we have

$$y_{m-1}^k - x^* - F_c'\left(x^k\right)^{-1} F_c\left(y_{m-1}^k\right) = o\left(\left(r^k\right)^{m+1}\right). \tag{21}$$

The following proves that holds for *m*, namely

$$r^{k+1} = y_m^{k+1} - x^* - F_c'\left(x^k\right)^{-1} F_c\left(y_{m-1}^k\right) = o\left(\left(r^k\right)^{m+2}\right). \tag{22}$$

Combine formula (20) and assumption (21), we can get

$$r_m^k = y_m^k - x^* = o\left(\left(r^k\right)^{m+1}\right), \tag{23}$$

$$F_c'\left(x^k\right)r^{k+1} = F_c'\left(x^k\right)r^k - F_c\left(y_m^k\right). \tag{24}$$

Using the expansion of $F_c'\left(x^k\right)$ and $F_c\left(y_m^k\right)$ at $x^*$, we obtain

$$F_c'\left(x^*\right)\left[I + o\left(r^k\right)\right]r^{k+1} = F_c'\left(x^*\right)\left[I + o\left(r^k\right)\right]r_m^k - F_c'\left(x^*\right)r_m^k - o\left(\left(r_m^k\right)^2\right). \tag{25}$$

From (24), then we got $r^{k+1} = o\left(\left(r^k\right)^{m+2}\right)$. The theorem is proved.

## 4. Numerical Results

In this section, we use numerical examples to examine the numerical effectiveness of smoothing modulus-based iterative methods from aspects of number of iteration steps (denoted by "IT"), elapsed CPU time in seconds (denoted by "CPU"), and norm of absolute residual vectors (denoted by "RES"). RES is defined as:

$$\text{RES} = abs\left(\left(Mz^k + q\right)^{\mathrm{T}}\left(z^k - m\left(z^k\right)\right)\right),$$

where $z^k$ is the $k$th approximate solution to the ICP.

In this paper, all calculations are run on a machine with a CPU of 1.8 Hz and a memory of 8G, and the programming language is MATLAB (2018b). We choose the $z^0 = (0, 0, \cdots, 0, 0)^{\mathrm{T}}$, large positive integer $c = 30$, $\xi_1 = \xi_2 = 1$, $\varepsilon = 10^{-6}$, $\alpha = \beta = 1$.

We will compare our smoothing modulus-based iterative method with accelerated modulus-based matrix splitting iteration methods, and its iteration coefficient is 1.6 [15]. The abbreviations of methods are listed in Table 1.

**Example 4.1.** Let $p$ is a positive integer, $n = p^2$, consider the implicit complementarity problem (1), when $M \in R^{n \times n}$ is a tridiagonal block matrix, $q$ is a vector

$$M = \begin{pmatrix} S & -I & 0 & \cdots & 0 \\ -I & S & \ddots & \ddots & \vdots \\ 0 & \ddots & \ddots & \ddots & 0 \\ \vdots & \ddots & \ddots & S & -I \\ 0 & \cdots & 0 & -I & S \end{pmatrix} \in R^{n \times n}, \quad q = \begin{pmatrix} -1 \\ 1 \\ \vdots \\ (-1)^{n-1} \\ (-1)^n \end{pmatrix} \in R^n$$

where $S = tridiag(-1, 4, -1) \in R^{p \times p}$ is a tridiagonal matrix, $I \in R^{p \times p}$ is an identity matrix, the point-to-point mapping $m(z) = \left(\sqrt{z_1}, \sqrt{z_2}, \cdots, \sqrt{z_n}\right)^{\mathrm{T}}$. The numerical results are listed in Table 2.

**Example 4.2.** Let $p$ is a positive integer, $n = p^2$, consider the implicit complementarity problem (1), when $M \in R^{n \times n}$ is a tridiagonal block matrix, $q$ is a vector

$$M = \begin{pmatrix} S & -0.5I & 0 & \cdots & 0 \\ -1.5I & S & \ddots & \ddots & \vdots \\ 0 & \ddots & \ddots & \ddots & 0 \\ \vdots & \ddots & \ddots & S & -0.5I \\ 0 & \cdots & 0 & -1.5I & S \end{pmatrix} \in R^{n \times n}, \quad q = \begin{pmatrix} -1 \\ 1 \\ \vdots \\ (-1)^{n-1} \\ (-1)^n \end{pmatrix} \in R^n$$

where $S = tridiag(-1.5, 4, -0.5) \in R^{p \times p}$ is a tridiagonal matrix, $I \in R^{p \times p}$ is an identity matrix, the point-to-point mapping $m(z) = \left(\arctan(z_1), \arctan(z_2), \cdots, \arctan(z_n)\right)^{\mathrm{T}}$. The numerical results are listed in Table 3.

From Table 2 and Table 3, for different problem of size $n$, we list the iteration steps, the CPU time and the residual norms with respect to AMSOR, SMN, MSMN and SM($m + 1$) methods. It can be seen that all methods converge quickly. Among these methods, SMN, MSMN and SM($m + 1$) methods require

**Table 1.** Abbreviations of methods.

| | |
|---|---|
| AMSOR | Accelerated Modulus-based Matrix Splitting Iteration Methods |
| SMN | Smoothing Modulus-based Newton Methods |
| MSMN | Modified Smoothing Modulus-based Newton Methods |
| SM$(m+1)$ $(m=3)$ | Smoothing Modulus-based $m+1$-step Iterative Methods |

**Table 2.** Numerical results of Example 4.1.

| Method | $n=400$ | | | $n=1600$ | | | $n=3600$ | | | $n=6400$ | | |
|---|---|---|---|---|---|---|---|---|---|---|---|---|
| | IT | CPU | RES | IT | CPU | RES | IT | CPU | RES | IT | CPU | RES |
| AMSOR | 5 | 0.038 | 3.41e−07 | 4 | 0.59 | 7.51e−07 | 3 | 3.25 | 5.50e−07 | 4 | 20.14 | 3.75e−07 |
| SMN | 5 | 0.021 | 2.21e−05 | 7 | 0.32 | 5.48e−05 | 8 | 1.98 | 1.03e−04 | 8 | 4.96 | 1.29e−04 |
| MSMN | 3 | 0.027 | 8.84e−16 | 3 | 0.21 | 2.45e−14 | 3 | 0.83 | 6.19e−14 | 3 | 2.88 | 1.10e−13 |
| SM$(m+1)$ | 2 | 0.017 | 2.73e−10 | 3 | 0.20 | 7.96e−14 | 3 | 0.78 | 7.65e−13 | 3 | 2.54 | 3.02e−09 |

**Table 3.** Numerical results of Example 4.2.

| Method | $n=400$ | | | $n=1600$ | | | $n=3600$ | | | $n=6400$ | | |
|---|---|---|---|---|---|---|---|---|---|---|---|---|
| | IT | CPU | RES | IT | CPU | RES | IT | CPU | RES | IT | CPU | RES |
| AMSOR | 5 | 0.040 | 4.64e−07 | 4 | 0.55 | 2.46e−07 | 3 | 3.54 | 1.18e−07 | 4 | 13.06 | 1.08e−07 |
| SMN | 2 | 0.025 | 2.81e−05 | 6 | 0.37 | 8.69e−05 | 6 | 1.43 | 1.64e−04 | 7 | 4.12 | 2.57e−04 |
| MSMN | 3 | 0.022 | 2.13e−15 | 3 | 0.18 | 3.09e−15 | 3 | 0.75 | 1.51e−14 | 3 | 2.41 | 3.19e−14 |
| SM$(m+1)$ | 2 | 0.021 | 3.00e−10 | 2 | 0.14 | 1.23e−09 | 2 | 0.69 | 3.50e−09 | 2 | 1.83 | 8.83e−09 |

less iteration steps and cost less computing time than AMSOR method, and SM$(m+1)$ is best.

**Example 4.3.** Let $p$ is a positive integer, $n=p^2$, consider the implicit complementarity problem (1), when $M \in R^{n \times n}$ is a tridiagonal block matrix, $q$ is a vector

$$M = \begin{pmatrix} S & -I & 0 & \cdots & 0 \\ -I & S & \ddots & \ddots & \vdots \\ 0 & \ddots & \ddots & \ddots & 0 \\ \vdots & \ddots & \ddots & S & -I \\ 0 & \cdots & 0 & -I & S \end{pmatrix} \in R^{n \times n}, \quad q = \begin{pmatrix} -1 \\ 1 \\ \vdots \\ (-1)^{n-1} \\ (-1)^n \end{pmatrix} \in R^n$$

where $S = tridiag(-1, 4, -1) \in R^{p \times p}$ is a tridiagonal matrix, $I \in R^{p \times p}$ is an identity matrix, the point-to-point mapping $m(z) = \left(z_1^3, z_1^3, \cdots, z_n^3\right)^T$. The numerical results are listed in Table 4.

**Example 4.4.** Let $p$ is a positive integer, $n=p^2$, consider the implicit complementarity problem (1), when $M \in R^{n \times n}$ is a tridiagonal block matrix, $q$ is a vector

**Table 4.** Numerical results of Example 4.3.

| Method | $n = 3025$ | | | $n = 24025$ | | | $n = 42025$ | | |
|--------|----|-----|-----|----|--------|----------|----|--------|----------|
| | IT | CPU | RES | IT | CPU | RES | IT | CPU | RES |
| SMN | 8 | 0.576 | 7.98e−04 | 10 | 26.472 | 2.806e−04 | 10 | 269.01 | 3.890e−04 |
| MSMN | 3 | 0.265 | 5.35e−14 | 2 | 8.690 | 2.265e−13 | 3 | 99.93 | 3.994e−10 |
| SM($m + 1$) | 2 | 0.192 | 4.42e−16 | 2 | 7.339 | 7.996e−14 | 2 | 70.03 | 2.706e−14 |

**Table 5.** Numerical results of Example 4.4.

| Method | $n = 3025$ | | | $n = 24025$ | | | $n = 42025$ | | |
|--------|----|-----|-----|----|--------|----------|----|--------|----------|
| | IT | CPU | RES | IT | CPU | RES | IT | CPU | RES |
| SMN | 6 | 0.432 | 1.397e−4 | 7 | 19.706 | 7.056e−4 | 7 | 240.62 | 1.132 e−4 |
| MSMN | 3 | 0.266 | 1.27e−14 | 3 | 9.798 | 1.143e−13 | 2 | 108.35 | 1.444e−14 |
| SM($m + 1$) | 2 | 0.214 | 9.65e−15 | 2 | 7.394 | 9.060e−14 | 2 | 77.316 | 1.688e−13 |

$$M = \begin{pmatrix} S & -0.5I & 0 & \cdots & 0 \\ -1.5I & S & \ddots & \ddots & \vdots \\ 0 & \ddots & \ddots & \ddots & 0 \\ \vdots & \ddots & \ddots & S & -0.5I \\ 0 & \cdots & 0 & -1.5I & S \end{pmatrix} \in R^{n \times n}, \quad q = \begin{pmatrix} -1 \\ 1 \\ \vdots \\ (-1)^{n-1} \\ (-1)^n \end{pmatrix} \in R^n$$

where $S = tridiag(-1.5, 4, -0.5) \in R^{p \times p}$ is a tridiagonal matrix, $I \in R^{p \times p}$ is an identity matrix, the point-to-point mapping $m(z) = (z_1^3, z_1^3, \cdots, z_n^3)^T$. The numerical results are listed in Table 5.

From Table 4 and Table 5, for the symmetric and asymmetric problem, we list the iteration steps, the CPU time and the residual norms with respect to SMN, MSMN and SM($m + 1$) methods respectively, among the three algorithm, SM($m + 1$) has the least number of iteration steps, costs the least CPU time, and holds better error accuracy.

## 5. Conclusion

In this paper, a class of smoothing modulus-based iteration method for solving implicit complementarity problems is proposed, the convergence of the algorithm is analyzed, and numerical experiments show the effectiveness of the method.

## Acknowledgements

This work was supported by Natural Science Foundation of China (11661027), the Guangxi Natural Science Foundation (2020GXNSFAA159143), and the Innovation Project of GUET Graduate Education (2021YCXS114).

## Conflicts of Interest

The authors declare no conflicts of interest regarding the publication of this paper.